\documentclass[12pt]{amsart}

\usepackage{a4wide}

\usepackage{amssymb}
\usepackage{epsfig}

\usepackage{caption}
\usepackage{subcaption}

\usepackage{enumerate}

\newcommand{\C}{\mathbb{C}}
\newcommand{\R}{\mathbb{R}}
\newcommand{\Z}{\mathbb{Z}}
\newcommand{\F}{\mathbb{F}}
\newcommand{\Q}{\mathbb{Q}}
\renewcommand{\H}{\mathbb{H}}
\renewcommand{\P}{\mathbb{P}}
\newcommand{\B}{\mathbb{B}}
\newcommand{\Sc}{\mathcal{S}}
\newcommand{\T}{\mathcal{T}}
\newcommand{\tr}{\mathrm{tr}}
\newcommand{\br}{\mathrm{br}}

\newcommand{\om}{\omega}
\newcommand{\omc}{\overline{\omega}}

\def\cqfd{\mbox{}\nolinebreak\hfill$\Box$\medbreak\par}
\newenvironment{pf}{\noindent\textbf{Proof:}}{\cqfd}

\newtheorem{thm}{Theorem}[section]
\newtheorem{prop}{Proposition}[section]
\newtheorem{dfn}{Definition}[section]

\theoremstyle{remark}
\newtheorem{rk}{Remark}[section]

\author{Martin Deraux}

\title[Non-arithmetic ball quotients]{Non-arithmetic ball quotients
  from a configuration of elliptic curves in an Abelian surface}

\address{Martin Deraux\\
Universit\'e Grenoble Alpes}
\email{martin.deraux@univ-grenoble-alpes.fr}

\date{June 16, 2017}

\begin{document}

\begin{abstract}
We construct some non-arithmetic ball quotients as branched covers of
a quotient of an Abelian surface by a finite group, and compare them
with lattices that previously appear in the literature. This gives an
alternative construction, which is independent of the computer, of
some lattices constructed by the author with Parker and Paupert.
\end{abstract}

\maketitle

\section{Introduction}

Recall that finite groups generated by complex reflections have been
classified by Shephard and Todd~\cite{shephardtodd}. Such a group $G$
comes with an isometric action on $\P^n$ (for the Fubini-Study
metric), and the quotient $X=G\setminus \P^n$ turns out to be a
weighted projective space. In fact the ring of invariant polynomials
in $n+1$ variables is generated by homogeneous polynomials
$f_{d_0},\dots,f_{d_{n}}$, where $f_{d_j}$ has degree $d_j$, hence the
quotient, which is given by the projective spectrum of the ring of
invariants, is the weighted projective plane $\P(d_0,\dots,d_{n})$.
Note that the weights can be computed from simple combinatorial data,
since the degrees satisfy $d_0\cdot d_1\cdots d_{n}=|G|$, and
$\sum(d_j-1)$ is equal to the number of reflections in the group
(see~\cite{springer}, for instance).

An analogous classification has been produced for affine
crystallographic complex reflection groups, see~\cite{popov}. The
basis for the classification of affine groups is the fact that the
group of automorphisms of affine space $\C^n$ is a semi-direct product
$V\rtimes GL(V)$, where $V$ is the vector space of translations in
$\C^n$, which allows us to reduce the classification to the problem of
classifying extensions of finite unitary reflection groups by a
lattice in $\C^n$. In particular, if $G$ is an affine crystallographic
group, the quotient $G\setminus \C^n$ is the quotient of a complex
torus by a finite group. Note that in general, $G$ is not the
semi-direct product of its linear part and its translation subgroup
(see p.57 in~\cite{popov}).

It was observed by Bernstein and Schwarzman~\cite{bernsteinSchwarzman}
that, at least in many cases, if $G$ is generated by complex
reflections, the quotient $G\setminus \C^n$ is again a weighted
projective plane. The heart of their proof is to construct suitable
$\Theta$-functions that play the role of the homogeneous invariant
polynomials in the Shephard-Todd-Chevalley theorem, which they managed
to do only when the linear part of $G$ is a real Coxeter group (in
that case, the weights of the weighted projective space are given by
the so-called \emph{exponents} of the corresponding Coxeter group).

Some quotients $G\setminus \C^2$, where $G$ is an affine
crystallographic complex reflection group whose linear part is not a
Coxeter group, were worked out by Kaneko, Tokunaga and
Yoshida~\cite{kanekotokunagayoshida}, building on the
Bernstein-Schwarzman result. The corresponding quotients turn out to
be explicit weighted projective planes, but their proof does not shed
much light on the general case.  Still, for a general affine
crystallographic complex reflection group $G$, it is believed that the
quotient $G\setminus\C^n$ should be a weighted projective space (see
p.~17 of~\cite{dolgachev}).

In this note, we investigate a particular affine crystallographic
complex reflection group $G$, whose linear part is the Shephard-Todd
group $G_{12}$, and whose subgroup of translations is given by the
lattice $\Lambda=(\Z\oplus i\sqrt{2}\Z)^2$. In other words, there is
an extension
\begin{equation}\label{eq:extension}
  1\rightarrow \Lambda\rightarrow G \rightarrow G_{12}\rightarrow 1,
\end{equation}
and one can think of the quotient $G\setminus \C^2$ as the quotient of
the Abelian surface $A=\C^2/\Lambda$ by the group $G_{12}$. For
concreteness, we mention that the group $G_{12}$ has order 48, it is a
central extension of the octahedral group, and it is also known to be
isomorphic to $GL(2,\F_3)$.

Our group $G$ is not a semi-direct product (equivalently the
sequence~\eqref{eq:extension} does not split), which characterizes it
uniquely up affine equivalence, according to~\cite{popov}. Note also
that the action of $G_{12}$ has no global fixed point in the Abelian
surface $A$, so our action is \emph{not} the same as the action given
by Birkenhake and Lange (see Theorem 13.4.5 in~\cite{birkenhakelange})

Since $G_{12}$ is not a Coxeter group, it is not in the list of groups
treated by Bernstein and Schwarzman, and it is not in the list of
groups treated by Kaneko, Tokunaga and Yoshida, so the structure of
the quotient seems to be unknown.

We will show that the quotient $X=G\setminus \C^2$ has two singular
points of type $\frac{1}{3}(1,2)$ and $\frac{1}{8}(1,3)$ respectively,
and that the map $A\rightarrow X$ ramifies with order 2 along a
(highly singular) rational curve in $X$, which is given by the image
of the union of all mirrors of complex reflections in the group $G$.
We refer to the branch locus as the discriminant curve, and denote it
by $M$. It does not contain any of the singular points of $X$, and the
curve $M$ has four singular points, two ordinary cusps, a point with
multiplicity four and another with multiplicity 6 (see
Figure~\ref{fig:models} for a schematic picture of the singularities
of $M$).

Assuming that $X$ is indeed a weighted projective plane, the list of
its singular points shows that it must be isomorphic to
$\P(1,3,8)$. The curve $M$ would then be an irreducible curve of
homogeneous degree 24, whose explicit equation remains elusive
(see~\cite{derauxklein} for the analogous equation in the case of
$\P(2,3,7)$ in relation to the Klein quartic).

We can rephrase the preceding paragraphs as follows.
\begin{thm} \label{thm:affine}
  The pair $(X,\frac{1}{2}M)$ is an orbifold which is uniformized by
  $\C^2$, and $\pi_1^{orb}(X,\frac{1}{2}M)$ is the affine
  crystallographic complex reflection group $G$.
\end{thm}
Theorem~\ref{thm:affine} parallels Proposition~2
of~\cite{derauxklein}, which says that $\P(2,3,7)$, with a specific
curve with weight $\frac{1}{2}=1-\frac{1}{2}$, gives an orbifold
uniformized by the positively curved complex space form $\P^2$. The
main result in~\cite{derauxklein} is obtained by changing $p=2$ to
higher integer values, i.e. changing the weight of the curve to be
$1-\frac{1}{p}$ (in fact, for most values of $p$, more subtle
modifications are needed).

It is then tempting to mimic the construction of~\cite{derauxklein},
and to consider the pairs $(X,(1-\frac{1}{p})M)$ for integer values
$p>2$ (in fact, it is convenient to allow also $p=\infty$). The basic
questions are the following.
\begin{enumerate}[(i)]
  \item When is the pair $(X,(1-\frac{1}{p})M)$ an orbifold? When it
    is an orbifold, is it modeled on a space form?
  \item When it is not an orbifold, is there a suitable model
    birational to it that is an orbifold? If so, is that orbifold
    modeled on a space form?
\end{enumerate}
The main goal of the present paper is to show that, even though the
answer to (i) (only for $p=2$) may seem disappointing, there is an
affirmative answer to question (ii) for some other well-chosen values
of $p$, namely $p=3,4,6$ or $\infty$. For these values, the universal
cover is the 2-dimensional complex space form of curvature $-1$, which
we denote by $\H^2$ (for basic facts on the complex hyperbolic plane
$\H^2$ and lattices in its isometry group, see section~\ref{sec:chg}).

The precise statements are somewhat technical (they will only be given in
section~\ref{sec:statement}), because on the one hand the birational
modifications are not that easy to describe, and on the other only a
proper open set turns out to be uniformized by $\H^2$. For now we
suggest that the reader keeps in mind that the statements below
roughly say that there is a indeed a complex hyperbolic uniformization
of suitable open sets in the pairs $(X,(1-\frac{1}{p})M)$, for
$p=3,4,6,\infty$.

For $p=3$, we will prove the following. 
\begin{thm} \label{thm:baby}
  The pair $(X,(1-\frac{1}{3})M)$ is a compactification of a ball
  quotient. More precisely, there is a lattice $\Gamma_3\subset
  PU(2,1)$ with one cusp, such that $X_0=\Gamma_3\setminus \H^2$ has
  1-point compactification isomorphic to $X$. Modulo this isomorphism,
  the quotient map $\H^2\rightarrow X_0$ branches with order $3$ along
  $M_0$, which is obtained from $M$ by removing its point with
  multiplicity 6.
\end{thm}
Note that the presence of a 1-dimensional branch locus for the
quotient map $\H^2\rightarrow X$ says that the lattice $\Gamma_3$
contains complex reflections, we will see later that it is actually
generated by complex reflections. Observe also that the point that needs
to be removed from $X$ in order to get a ball quotient is
characterized by the fact that it is the only point where the pair
$(X,(1-\frac{1}{3})M)$ is not log-terminal, see
section~\ref{sec:statement}.

The orbifold structure on $X=A/G$ does not lift to an orbifold
structure on the Abelian surface $A$, since the corresponding weights
on the preimage of the discriminant curve would have to be equal to
$3/2$, which is not an integer.

There are statements analogous to Theorem~\ref{thm:baby} for $p=4$,
$6$ or $p=\infty$, but then the pair $(X,(1-\frac{1}{p})M)$ is
actually not an orbifold, and one needs to perform a suitable
birational modification before it becomes one. After suitable
modification, for each case $p=4$, $6$ or $p=\infty$, one gets an
orbifold which is uniformized by $\H^2$, with orbifold fundamental
group given by a non-cocompact lattice $\Gamma_p$.

For now we simply give a rough statement.
\begin{thm} \label{thm:rough}
  There are a lattices $\Gamma_p\subset PU(2,1)$, $p=4,6,\infty$ such
  that $\Gamma_p\setminus \H^2$ has a compactification birational
  to $X$. The groups $\Gamma_p$ have one cusp for $p=4,\infty$, two
  cusps for $p=6$.
\end{thm}
The explicit birational transformation that yields the corresponding
compactification will be given later in the paper (see
section~\ref{sec:statement}, Theorem~\ref{thm:main} in particular).

For the groups that appear in Theorem~\ref{thm:rough}, the weight of
the orbifold structure along (the strict transform of) the
discriminant curve is even, so the orbifold structure on (the suitable
surface birational to) $X=A/G$ lifts to an orbifold structure on (a
suitable blow up of) the Abelian surface $A$, with multiplicity
$2=4/2$, $3=6/2$ or $\infty$ respectively at a generic point of the
union of mirrors of reflections of $G$. The fact that the orbifold
structure lifts to $A$ only when $p$ is even has a similar incarnation
in Deligne-Mostow theory, when passing from the Picard integrality
condition INT to the condition $\Sigma$-INT
(see~\cite{delignemostow},~\cite{mostowihes},~\cite{delignemostowbook}).

For $p=3,4$ and $6$, the lattices $\Gamma_p$ turn out to be conjugate
to lattices constructed by the author in joint work with Parker
and Paupert, see~\cite{dpp1} and~\cite{thealgo}, namely the
groups $\Sc(p,\sigma_1)$, generated by a complex reflection $R_1$ of
angle $2\pi/p$, and a regular elliptic element $J$ of order 3 such
that
$
\tr(R_1J)=-1+i\sqrt{2}.
$ 
For basic notation on these groups, see section~\ref{sec:sporadic}.

It was proved in~\cite{dpp1} and~\cite{thealgo} that $\Sc(p,\sigma_1)$
is discrete if and only if $p=3,4,6$, and in those cases it is a
non-cocompact lattice. It has one cusp for $p=3,4$, two cusps for
$p=6$. Note also that the three groups can be checked to be generated
by complex reflections, namely by $R_1$, $R_2=JR_1J^{-1}$ and
$R_3=J^{-1}R_1J$ (see section~\ref{sec:sporadic}).

We will prove the following.
\begin{thm}
  For every $p=3,4,6$, the group $\Gamma_p$ is conjugate in $PU(2,1)$
  to the group $\Sc(p,\sigma_1)$.
\end{thm}
In particular, because of the analysis in~\cite{thealgo}, we know that
the $\Gamma_p$, $p=3,4,6$ are non-arithmetic lattices. The group
corresponding to $p=\infty$ does not appear in~\cite{dpp1}, but it is
in a sense less interesting since it turns out to be arithmetic.

Complex hyperbolic lattices have been previously constructed
from configurations of elliptic curves on an Abelian surface. One
important construction was worked out by Livne, see~\cite{livne} (and
also~\cite{delignemostowbook}), from a point of view that is fairly
different from ours. Another construction, closer in spirit to
the results in this paper, appears in~\cite{hirzabelian} (see
also~\cite{stover},~\cite{dicerbostover},~\cite{roulleauurzua} for
recent developments).

Just as in~\cite{derauxklein}, the results of this paper give an
alternative construction of certain non-arithmetic ball quotients,
whose existence was known so far only by giving explicit matrix
generators and constructing a fundamental domain for their action
(see~\cite{dpp1} and~\cite{thealgo}). 

The analysis in~\cite{derauxklein} shows that some of the
non-arithmetic lattices in~\cite{thealgo}, even though they are not
commensurable to Deligne-Mostow lattices
(see~\cite{delignemostow},~\cite{mostowihes}), are commensurable to
Couwenberg-Heckman-Looijenga lattices (see~\cite{chl}, which was
inspired in part by~\cite{bhh}). For brevity, we refer to these two
classes of lattices as DM and CHL lattices, respectively (note that DM
lattices are special cases of CHL lattices). In fact, an analysis
similar to the one in~\cite{derauxklein} shows the following (for
notation of Sporadic and Thompson triangle groups, see
section~\ref{sec:sporadic}).
\begin{thm}\label{thm:chl}
\begin{enumerate}
\item The group $\mathcal{S}(2,\sigma_{10})$ is isomorphic to the
  Shephard-Todd group $G_{23}$. The lattices
  $\mathcal{S}(p,\sigma_{10})$, $p=3,4,5,10$ are conjugate to the
  corresponding CHL lattices of type $H_3$.
\item The group $\mathcal{S}(2,\overline{\sigma}_{4})$ is isomorphic
  to a subgroup of index two in the Shephard-Todd group $G_{24}$, both
  groups having isomorphic projectivizations of order 168. The
  lattices $\mathcal{S}(p,\overline{\sigma_{4}})$, $p=3,4,5,6,8,12$ are
  conjugate to the corresponding CHL lattices.
\item The group $\mathcal{T}(2,{\bf S_2})$ is isomorphic to a subgroup
  of index two in the Shephard-Todd group $G_{27}$, both having
  isomorphic projectivization of order 360. The lattices
  $\mathcal{T}(p,{\bf S_2})$, $p=3,4,5$ are conjugate to the
  corresponding CHL lattices.
\end{enumerate}
\end{thm}
The three families of lattices in Theorem~\ref{thm:chl}, together
with Deligne-Mostow lattices, exhaust the list of CHL lattices in
$PU(2,1)$ (the other ones constructed in~\cite{chl} are in $PU(n,1)$
for $n>2$).

In particular, we have the following.
\begin{thm}
  The lattices $\Sc(p,\sigma_1)$, $p=3,4,6$ are not commensurable to
  any CHL lattice (and in particular not to any DM lattice either).
\end{thm}

Some lattices in~\cite{thealgo} are still not treated by the methods
in~\cite{derauxklein} nor of the present paper, for instance the
sporadic lattices $\mathcal{S}(p,\sigma_5)$. Indeed, in the family of
$\sigma_5$ groups, there seems to be no finite nor any
crystallographic group.\\

\noindent\textbf{Acknowledgements:} I wish to thank St\'ephane Druel
for many stimulating discussions related to the results in this paper,
as well as Xavier Roulleau for several comments that significantly
helped improve the exposition.

\section{Basic complex hyperbolic geometry} \label{sec:chg}

Recall that the complex hyperbolic plane $\H^2$ is the only complete,
simply connected K\"ahler surface of holomorphic sectional curvature
$-1$. It is biholomorphic to the unit ball $\B^2\subset \C^2$, and we
equip it with the only metric that is invariant under the group of
biholomorphisms of $\B^2$ (normalized so that the holomorphic
sectional curvature is $-1$). In terms of Riemannian symmetric spaces,
$\H^2$ is the non-compact dual of $\P^2$. We summarize a few basic
facts that we will use in this paper (see~\cite{goldman} for much more
information).

Working in homogeneous coordinates for $\P^2$ and seeing
$\B^2\subset\C^2\subset \P^2$ as sitting in an affine chart of the
complex projective plane, one can see biholomorphisms of $\B^2$ as
induced by linear transformations of $\C^3$ that preserve a Hermitian
form of signature $(2,1)$, say $\langle Z,W\rangle = -Z_0\bar
W_0+Z_1\bar W_1+Z_2\bar W_2$. The unit ball is then identified with
the set of negative complex lines in $\C^3$, i.e. lines spanned by a
vector $V$ with $\langle V,V\rangle<0$. This description gives an
isomorphism $Bihol(\B^2)\simeq PU(2,1)$, which produces almost all
isometries of $\H^2$ (the full group of isometries is generated by
$PU(2,1)$ and the single isometry given by complex conjugation).

We will use the classification of (non-trivial) isometries of
$\H^2$ into elliptic, parabolic and loxodromic elements
(see~\cite{chengreenberg} for instance). Elliptic isometries are
characterized by the fact that they fix at least one point in
$\H^2$. Parabolic elements have unique fixed point at infinity,
i.e. in $\partial_{\infty}\H^2\simeq S^3$. Loxodromic elements have
precisely two fixed points at infinity.

Elliptic isometries whose matrix representatives have distinct
eigenvalues are called regular elliptic isometries. Among non-regular
elliptic isometries, an important class is given by complex
reflections, that fix pointwise the intersection with $\B^2$ of an
affine complex line in $\C^2$. These are characterized in terms of
their matrix representative in $U(2,1)$ by the fact that they have a
double eigenvalue, and that the simple eigenvalue eigenspace is
spanned by a vector with positive square norm.

A lattice $\Gamma\subset PU(2,1)$ is a discrete subgroup such that
$\Gamma\setminus PU(2,1)$ has finite Haar measure. Equivalently, the
quotient $\Gamma\setminus \H^2$ has finite volume for the Riemannian
volume form on $\H^2$. $\Gamma$ is called co-compact (or uniform) if
the quotient $\Gamma\setminus \H^2$ is compact. If it is not, there
are finitely many conjugacy classes of maximal parabolic subgroups in
$\Gamma$, and the quotient decomposes as a disjoint union of a compact
part and finitely many cusps (a cusp is the quotient of a sufficiently
small horoball centered at the fixed point of one of the parabolic
subgroups). We say $\Gamma$ has $n$ cusps if the quotient has $n$
cusps, equivalently if there are $n$ conjugacy classes of maximal
parabolic subgroups in $\Gamma$.

\section{Complex hyperbolic lattice triangle groups} \label{sec:sporadic}

\subsection{Sporadic triangle groups}

In this section we briefly review some of the basic facts and notation
in~\cite{thealgo} (and also previous papers cited there). In the
following statement, we write $\omega=(-1+i\sqrt{3})/2$.
\begin{prop}\label{prop:definegroups}
  Let $p\in \R^*$, $u=e^{2\pi i/3p}$, $\tau,\tau'\in\C$. Up to conjugacy in $SL(3,\C)$,
  there is a unique pair $(R_1,J)$ of matrices such that
  \begin{itemize}
  \item $R_1$ has eigenvalues $u^2,\bar u,\bar u$;
  \item $J$ has eigenvalues $1,\omega,\bar\omega$;
  \item $\tr(R_1J)=\tau$ and $\tr(R_1J^{-1})=\tau'$.
  \end{itemize}
  The group generated by $R_1$ and $J$ preserves a non-zero Hermitian
  form if and only if $\tau'=-u\overline{\tau}$, and in that case the form is
  unique up to scaling.
\end{prop}
Choosing the basis of $\C^3$ given by $e_1,e_2=Je_1,e3=J^{-1}e_1$, we can write
\begin{equation}\label{eq:R1J}
R_1 = \left(
  \begin{matrix}
    u^2 & \tau & \tau'\\
    0   &\bar u&  0\\
    0   & 0    &\bar u
  \end{matrix}\right), 
\quad J = \left(
  \begin{matrix}
    0 & 0 & 1\\
    1 & 0 & 0\\
    0 & 1 & 0
  \end{matrix}\right). 
\end{equation}
In the Hermitian case, i.e. when $\tau'=-u\overline{\tau}$, and assuming moreover
that $u^3\neq 1$, the invariant Hermitian form is given (up to a
nonzero scalar) by
\begin{equation}\label{eq:H}
  \left(
  \begin{matrix}
    \alpha & \beta & \bar\beta\\
    \bar \beta&\alpha&\beta\\
    \beta&\bar\beta&\alpha
  \end{matrix}\right),
\end{equation}
where $\alpha=2-u^3-\bar u^3$, $\beta=(\bar u^2-u)\tau$. Note that
this matrix tends to 0 when $p\rightarrow +\infty$ (recall $u=e^{2\pi
  i/3p}$), but after rescaling it by $1/\sqrt{2-u^3-\bar u^3}$, it
converges to
  $$ 
  \left(
  \begin{matrix}
    0         & -i\tau    & i\bar\tau\\
    i\bar\tau &    0      & -i\tau\\
    -i\tau    & i\bar\tau & 0
  \end{matrix}\right),
  $$ 
which gives the invariant Hermitian form when $u^3=1$.
\begin{dfn}\label{dfn:sporadic}
  We denote by $\Sc(p,\tau)$ the group generated by $R_1$ and $J$ as
  in~\eqref{eq:R1J}, where $\tau'=-u\overline{\tau}$, and refer to it as a
  \emph{sporadic triangle group} with trace parameter $\tau$.
\end{dfn}
In such a sporadic triangle group, it is natural to consider
$$
R_2=JR_1J^{-1},\quad R_3=J^{-1}R_1J.
$$
The groups are constructed so that $R_1J$ has finite order (its order
is actually independent of $p$). When that order is not a multiple
of 3, the group generated by $R_1,R_2$ and $R_3$ is actually the same
as the group generated by $R_1$ and $J$ (in particular, in those
cases, $\Sc(p,\tau)$ is generated by complex reflections).

The main groups of interest in this paper will be the groups
$\Sc(p,\sigma_1)$, for suitable integer values of $p$, and
$$
\sigma_1=-1+i\sqrt{2}.
$$ 

The following is easily obtained using~\eqref{eq:H} and the above
discussion.
\begin{prop}
For $p\geq 2$ an integer, the Hermitian form preserved by
$\Sc(p,\sigma_1)$ is definite if and only if $p=2$, and that it has
signature $(2,1)$ for all $p>2$. For every $p$, $R_1J$ has order 8,
and the group $\Sc(p,\sigma_1)$ is generated by $R_1$, $R_2$ and
$R_3$.
\end{prop}
Note that $R_1J$ having order 8 is easily seen to imply
that $J=R_1R_2R_3R_1R_2R_3R_1R_2$.

\subsection{Thompson triangle groups}

The groups $\T(p,{\bf T})$ are analogs of the sporadic groups, that
were constructed in James Thompson's
Ph.D. thesis~\cite{thompson}. The are generated by three complex
reflections $R_1$, $R_2$ and $R_3$, that have the same rotation
angles, but are not cyclically conjugated by any element of $J$ order
3. Here ${\bf T}=(\rho,\sigma,\tau)$ is a triple of complex numbers that
generalize the trace parameter of sporadic triangle groups, related to
traces of $R_jR_k$. Since they are not central to this paper, we omit
the detailed description of these groups and simply refer
to~\cite{thealgo}.

The Thompson triangle groups that appear in Theorem~\ref{thm:chl} are
the groups $\T(p,{\bf S_2})$, with trace parameter triple ${\bf
  S_2}=(1+\omega\frac{1+\sqrt{5}}{2},1,1)$. 

We will also use the description of $\Sc(p,\sigma_1)$ as 3,3,4;6
triangle groups, in other words, in the terminology of~\cite{thealgo},
as $\mathcal{T}(p,{\bf E_1})$, where ${\bf E_1}=(i\sqrt{2},1,1)$.

Recall that the integers in 3,3,4;6 stand for specific braid lengths
$\br(a,b)$, namely $\br(R_2,R_3)=3$, $\br(R_3,R_1)=3$,
$\br(R_1,R_2)=4$, $\br(R_1,R_3^{-1}R_2R_3)=6$, and $\br(a,b)=k$ means
$$
(ab)^{k/2}=(ba)^{k/2},
$$
but $(ab)^{n/2} \neq (ba)^{n/2}$ for every $n<k$.

In other words, $\mathcal{T}(p,{\bf E_1})$ is a group generated by
three reflections $R_1$, $R_2$, $R_3$ of the same order $p$, such that
$(R_1R_2)^2=(R_2R_1)^2$, $R_2R_3R_2=R_3R_2R_3$, $R_3R_1R_3=R_1R_3R_1$,
$(R_1\cdot R_3^{-1}R_2R_3)^3=(R_3^{-1}R_2R_3\cdot R_1)^3$.

The fact that $\Sc(p,\sigma_1)$ is conjugate to $\T(p,{\bf E_1})$
follows from a change of generators along the same lines as
in~\cite{kamiyaparkerthompson} (for details, see section~7.1
of~\cite{thealgo}). Explicitly, if $M_1$, $M_2=JM_1J^{-1}$,
$M_3=J^{-1}M_1J$ denote standard generators for $\Sc(p,\sigma_1)$,
then the matrices
$$
R_1=(M_3M_1M_2M_1^{-1})M_3(M_3M_1M_2M_1^{-1})^{-1}, R_2=(M_3M_1)M_2(M_3M_1)^{-1}, R_3=M_1
$$ 
give another generating set, which exhibits an isomorphism with
$\T(p,{\bf E_1})$.

\section{The affine crystallographic reflection group} \label{sec:acrg}

We start by describing the relevant affine crystallographic group. One
way to write it is to use the matrices given in~\cite{dpp1}, which
would give slightly complicated matrices, and then to diagonalize the
corresponding Hermitian form by a suitable coordinate change. Here we
only give the matrices in a nice basis. For computational
convenience, rather than choosing the generators to have determinant
one as we did in~\cite{thealgo}, we adjust the repeated eigenvalue to
be equal to 1 (this amounts to multiplying the generators by a
suitable root of unity).

\begin{dfn}
Let $G$ be the group generated by the matrices $R_1$, $R_2$ and $R_3$ given below
{
\tiny
$$
R_1=\left(
\begin{matrix}
1 & 0 & 0\\
0 & 1 & 0\\
0 & 1-i\sqrt{2} & -1
\end{matrix}\right), \quad
R_2=\left(
\begin{matrix}
1 & 0 & 0\\
0 & -1+i\sqrt{2} & 2\\
0 & 1+i\sqrt{2} & 1-i\sqrt{2}
\end{matrix}\right), \quad
R_3=\left(
\begin{matrix}
1 & 0 & 0\\
\frac{1+i\sqrt{2}}{2} & 1 & -1-i\sqrt{2}\\
1 & 0 & -1
\end{matrix}\right), \quad
$$  
}
\end{dfn}
First observe that $G$ can be thought of as a subgroup of the
semi-direct product $\C^2\rtimes U(2)$. To see this, we write
$(z_0,z_1,z_2)$ for the coordinates in $\C^3$, and denote by
$\pi:\C^3\rightarrow\C^2$ the projection onto the last two
coordinates. Note that the group $G$ clearly preserves every
hyperplane $z_0=\lambda$, $\lambda\in\C$. We will study the affine
action of $G$ on $\C^2$ given by
$$
B\cdot(z_1,z_2)=\pi(B(1,z_1,z_2)).
$$ 
Concretely, we think of the linear part of $B$ as being given by
the lower right 2$\times$2 block of the 3$\times$3 matrix $B$, and the
translation part by the lower left 2$\times$1 block. We denote by
$\psi:G\rightarrow GL_2(\C)$ the corresponding homomorphism. Note that
the image of $\psi$ preserves a positive definite Hermitian form,
namely $\langle z,w\rangle = w^*Hz$ where
$$
H = \left(\begin{matrix}
 1 & \frac{-1-i\sqrt{2}}{2}\\
  \frac{-1+i\sqrt{2}}{2} & 1
\end{matrix}\right).
$$ 
The unitary group $U(H)$ is isomorphic to $U(2)$ since the matrix
$H$ has eigenvalues $2\pm \sqrt{3}$, which are both positive. One
checks that the matrices $\psi(R_j)$ are complex reflections of order
2 (i.e. each has eigenvalues 1 and -1), so $G$ is an affine group
generated by complex reflections. Next, we show that this group is
crystallographic, i.e. it is discrete, and the quotient of $\C^2$ by
its action is cocompact.  This follows from
Propositions~\ref{prop:finite} and~\ref{prop:crystal} below.

\begin{prop} \label{prop:finite}
  The linear part $\psi(G)$ of $G$ is isomorphic to the Shephard-Todd
  group $G_{12}$, which is isomorphic to $GL(2,\F_3)$.
\end{prop}

\begin{pf}
  One easily checks that the three matrices $A_1=\psi(R_1)$,
  $A_2=\psi(R_2)$ and $A_3=\psi(R_3)$ generate a group of order
  48, and that they satisfy the relations in a presentation for
  $G_{12}$, see~\cite{shephardtodd}, namely
  \begin{equation}
    A_1^2=A_2^2=A_3^2=(A_1A_2)^4=(A_2A_3)^3=(A_3A_1)^3=Id,
  \end{equation}
  and the element $(A_1A_2)^2$ is central of order 2.
\end{pf}

We refer to the matrix
$$
T_v=\left(
\begin{matrix}
1 & 0 & 0\\
v_1 & 1 & 0\\
v_2 & 0 & 1
\end{matrix}\right)
$$ 
as a translation with vector $v=(v_1,v_2)$. Let $K$ denote the
kernel of $\psi$, and let $T_{\Lambda}$ denote the group of
translations $T_v$, where $v\in\Lambda$ is a lattice vector.
\begin{prop} \label{prop:crystal}
  The groups $K$ and $T_{\Lambda}$ are equal.
\end{prop}

\begin{pf}
  One verifies that
  $$
  (R_3R_1R_2R_1)^2R_3R_2=T_{(i\sqrt{2},1)},\quad (R_3R_2R_1R_2)^2R_3R_1=T_{(0,1)}, 
  $$
  $$
  R_2[(R_2R_1)^2,R_3]R_2=T_{(-1,-1)},\quad (R_2R_1R_3R_1)^2R_1R_2R_3R_1=T_{(-1,i\sqrt{2})}.
  $$
  These four translations generate $T_{\Lambda}$, so we have $T_{\Lambda}\subset K$.

  In order to show the other inclusion, we observe that 
  \begin{equation}
    R_1^2=R_2^2=R_3^2=(R_1R_2)^4=(R_2R_3)^3=(R_3R_1)^3=Id,
  \end{equation}
  and the element $Z=(R_2R_1)^2$ commutes with $R_1$ and $R_2$, and it
  has order 2. The commutator $ZR_3Z^{-1}R_3^{-1}=(ZR_3)^2$ is given
  by the translation $T_v$, where $v=(-1-i\sqrt{2},-2)$. Once again,
  using the Shephard-Todd presentation for $G_{12}$, we get that
  $G/T_{\Lambda}$ is a quotient of $G_{12}$, but since
  $T_{\Lambda}\subset K$ and $G/K$ has order 48, both quotients of $G$
  must have order 48.
\end{pf}

In what follows, we denote by $F$ the finite group $G_{12}$. We denote
by $A$ the Abelian variety $\C^2/\Lambda$, and by $X$ the quotient of
$A$ by the action of $F=G/T_{\Lambda}$.

The following two propositions follow from painful (but not
particularly difficult) computation and bookkeeping.
\begin{prop}
  The group $F$ contains precisely 12 reflections, all of order 2,
  whose fixed point sets are elliptic curves in $A$. The group $F$
  acts transitively on the set of these 12 elliptic curves.
\end{prop}
For completeness, we list equations for these elliptic curves in
Table~\ref{tab:mirrorequations}.
\begin{table}[htbp]
  \begin{eqnarray*}
    1 & z_2=\frac{1-i\sqrt{2}}{2}z_1\\
    2 & z_2=\frac{2-i\sqrt{2}}{2}z_1\\
    3 & z_2=\frac{1}{2}\\
    121 & z_2=\frac{-i\sqrt{2}}{2}z_1\\
    131 & z_2=(1-i\sqrt{2})z_1-\frac{1}{2}\\ 
    212 & z_1=0\\
    232 & z_1=\frac{1+2i\sqrt{2}}{3}z_2+\frac{1-i\sqrt{2}}{6}\\
    32121 & z_1=\frac{1+i\sqrt{2}}{2}z_2+\frac{1+i\sqrt{2}}{4}\\
    23121 & z_1=z_2+\frac{1+i\sqrt{2}}{2}\\
    21321 & z_1=\frac{i\sqrt{2}}{2}z_2+\frac{2+i\sqrt{2}}{4}\\
    12321 & z_1=\frac{2+i\sqrt{2}}{2}z_2+\frac{i\sqrt{2}}{4}\\
   21231 & z_1=(1+i\sqrt{2})z_2+\frac{1+i\sqrt{2}}{2}\\
  \end{eqnarray*}
  \caption{Equations in $\C^2$ of (representatives of) the 12 mirrors
    of reflections in $A=\C^2/\Lambda$.}\label{tab:mirrorequations}
\end{table}
We denote by $\widetilde{M}$ the union of the mirrors in $A$, and by
$M$ its image in $X$. By transitivity of the action, $M$ is an
irreducible curve in $X$. 
\begin{prop} \label{prop:singsmooth}
  The action of $F$ on $A$ has precisely two orbits of fixed points in
  $A\setminus \widetilde{M}$, one with isotropy group of order 3, the
  other with isotropy group of order 8, as in
  Table~\ref{tab:sing}. The isotropy groups of points in
  $\widetilde{M}$ are all generated by complex reflections, the
  generic point having isotropy of order 2. The points in
  $\widetilde{M}$ with isotropy of order larger than 2 consist of two
  orbits of points with isotropy group of order 6, one orbit of points
  with isotropy of order 8, and one orbit of points with istropy of
  order 12, see Table~\ref{tab:smooth}.
\end{prop}
\begin{table}[htbp]
\begin{tabular}{|c|c|c|c|}
\hline
  Group                       & Order & eigenvalues          & coords\\
\hline
  $\langle R_1R_2R_3 \rangle$ & 8     & $\zeta_8,\zeta_8^3$  & $(\frac{1}{2},\frac{1+i\sqrt{2}}{2})$\\
  $\langle R_1R_3 \rangle$    & 3     & $\om,\omc$           & $(\frac{1+i\sqrt{2}}{3},\frac{1+2i\sqrt{2}}{6})$\\
\hline
\end{tabular}
\caption{Representatives of the orbits of points with non-reflection
  stabilizer (these produce singular points of the
  quotient).}\label{tab:sing}
\end{table}
\begin{table}[htbp]
\begin{tabular}{|c|c|c|c|c|c|}
\hline
  Generators          & Order  & ST-group   & Sing. of $M$  & notation        & coords\\
\hline
  $R_1,R_3$           & 6      & $G(3,3,2)$ & $z_1^3=z_2^2$            & $p_{13}$        & $(\frac{1+i\sqrt{2}}{3},\frac{1}{2})$\\
  $R_2,R_3$           & 6      & $G(3,3,2)$ & $z_1^3=z_2^2$            & $p_{23}$        & $(\frac{-2-i\sqrt{2}}{6},\frac{1}{2})$ \\
  $R_1,R_2$           & 8      & $G(2,1,2)$ & $z_1^4=z_2^2$            & $p_{12}$        & $(0,0)$\\
  $R_1,R_3(R_2R_1)^2$ & 12     & $G(6,6,2)$ & $z_1^6=z_2^2$            & $p_{13(21)^2}$  & $(0,\frac{1}{2})$\\ 
\hline
\end{tabular}
\caption{Representatives of the orbits of points whose stabilizer is a
  reflection group (these produce smooth points of the
  quotient).}\label{tab:smooth}
\end{table}

The results in Proposition~\ref{prop:singsmooth} follow by explicit
calculations. For a definition of the groups $G(m,p,n)$,
see~\cite{shephardtodd}, and also \S1 of~\cite{kanekotokunagayoshida},
for instance. The local analytic structure of the branch locus of the
quotient (fourth column in Table~\ref{tab:smooth}) can be obtained by
computing explicit invariant polynomials for the group, the results
are tabulated in~\cite{bannai}.

It follows from Proposition~\ref{prop:singsmooth} that $X$ has exactly
two singular points. Let $V$ denote the subset of $A$ of points with
trivial isotropy for the $F$-action, and let $U$ denote its image in
$X$.
\begin{prop}\label{prop:euler}
  We have $\chi(V)=48$, hence $\chi(U)=1$.
\end{prop}
\begin{pf}
  There are 48/3=16 points above the isolated singularity of order 3,
  48/8=6 points above the isolated singularity of order 8. There are
  $2\cdot(48/6)+(48/8)+(48/12)=16+6+4=26$ points with reflection
  isotropy of order $>2$. This gives 48 points.

  There are also 12 mirrors, each being an elliptic curve and
  containing 8 special points. The Euler characteristic of the
  generic stratum of each mirror is then $-8$, so we get
  $$
  0=\chi(A)=48+12\cdot(-8)+\chi(V),
  $$ 
  hence $\chi(V)=48$, and $\chi(U)=\chi(V)/48=1$, since $F$ has
  order 48.
\end{pf}

We will also need to study the stabilizer of a mirror of reflections.
\begin{prop}
  Each mirror in the group contains precisely 8 points with special
  isotropy (i.e. stabilizer of order strictly larger than 2).  The
  curve $M$ is a $\P^1$ with two pairs of points identified, and the
  map from each irreducible component of $\widetilde{M}$ to $M$ is a
  branched cover of degree 2.
\end{prop}
\begin{pf}
  In the coordinates we used above, the mirror of $R_2R_1R_2$
  corresponds to the elliptic curve $z_1=0$. The intersections with
  the other mirrors can be computed explicitly from the equations in
  Table~\ref{tab:mirrorequations}, they are listed in Table~\ref{tab:mirror}.
  \begin{table}[htbp]
    \begin{tabular}{|c|c|}
      \hline Mirrors & $z_2$\\ \hline 1,2,121,212 & $0$,
      $-i\sqrt{2}/2$\\ 212,232,12321 &
      $\pm(1+i\sqrt{2})/6$\\ 212,32121,21231 &
      $\pm(1+2i\sqrt{2})/6$\\ 1,3,131,212,32121,21231 &
      $1/2$\\ 1,212,232,23121,21321,12321 & $(1+i\sqrt{2})/2$\\ \hline
    \end{tabular}
    \caption{Special points on the mirror of $R_2R_1R_2$. We list the
      corresponding reflections whose mirrors meet at that point.}\label{tab:mirror}
  \end{table}
  One verifies that the only reflection that stabilizes the mirror of
  $R_2R_1R_2$ is $R_1$ (note that $R_1$ commutes with $R_2R_1R_2$,
  since $R_1R_2$ has order 4). Now $R_1$ acts on $z_1=0$ by
  $z_2\mapsto -z_2$. Among the points listed in
  Table~\ref{tab:mirror}, the two points $\pm(1+i\sqrt{2})/6$ get
  identified, and so do the points $\pm(1+2i\sqrt{2})/6$. The other
  four points are fixed by the action of $R_1$.
\end{pf}

\section{Statement of the main result}\label{sec:statement}

Recall that $X$ denotes the quotient $A/F$, where $A$ is the Abelian
variety $\C^2/\Lambda$, and $F$ is a specific group of order $48$,
isomorphic to the Shephard-Todd group $G_{12}$. As above, we denote by
$M\subset X$ the curve which is the image of the set of mirrors in $A$
of reflections of $F$.

We denote by $p_{12}$ the image in $X$ of the fixed point of $R_1R_2$,
etc (see Table~\ref{tab:smooth}). As in~\cite{derauxklein}, in order
to produce orbifolds uniformized by the ball, we will need to perform
suitable blow-ups on $X$.

The curve $M$ has a local analytic equation of the form
$(z_1^3-z_2)(z_1^3+z_2)=0$ near $p_{13(21)^2}$ (see
Table~\ref{tab:smooth}), so locally there are two tangent
components. The space $\widehat{Y}$ is obtained from $X$ by blowing up
$p_{13(21)^2}$ three times (the first blow-up preserves the tangency,
the second makes the intersection transverse, the third makes the two
local components disjoint). The exceptional locus of $\pi:\widehat{Y}\rightarrow
X$ is a chain of projective lines with self-intersections $-1,-2,-2$.
\begin{dfn}\label{def:Y}
  The space $Y$ is obtained from $\widehat{Y}$ by contracting the two
  $-2$ curves in the exceptional locus of $\pi:\widehat{Y}\rightarrow
  X$. We denote by $\gamma:\widehat{Y}\rightarrow Y$ the contraction,
  and $\varphi:Y\rightarrow X$ the corresponding birational
  transformation. We denote the exceptional locus of $\varphi$ by $E$.
\end{dfn}

Similarly, the space $\widehat{Z}$ is obtained from $X$ by blowing up
both points $p_{13(21)^2}$ and $p_{12}$. Near the first one, the
modification is the same as in the construction of $Y$. Near $p_{12}$,
the curve $M$ has a local equation of the form
$(z_1^2-z_2)(z_1^2+z_2)=0$ (see
Table~\ref{tab:smooth} again). At that point, we perform two successive
blow-ups (the first one makes the two tangent local components
transverse, the second makes them disjoint), which produces a chain of
two projective lines with self-intersection $-1,-2$.
\begin{dfn}\label{def:Z}
  The space $Z$ is obtained from $\widehat{Z}$ by contracting the two
  $-2$ curves above $p_{13(21)^2}$, and the $(-2)$-curve above
  $p_{12}$.  With a slight abuse of notation, we still denote by
  $\gamma:\widehat{Z}\rightarrow Z$ the contraction, and
  $\varphi:Z\rightarrow X$ the corresponding birational
  transformation. We denote the exceptional lines by $E$ and $F$,
  above $p_{13(21)^2}$ and $p_{12}$ respectively.
\end{dfn}

Finally, we consider $\widehat{W}$, which is obtained from $X$ by
blowing up all points $p_{13(21)^2}$ and $p_{12}$, $p_{13}$ and
$p_{23}$. Near each point $p_{13}$ and $p_{23}$, we need to perform
three successive blow-ups, producing a chain of three $\P^1$ with
self-intersections $-2,-1,-3$ (see~\cite{derauxklein}).
\begin{dfn}\label{def:W}
  The space $W$ is obtained from $\widehat{W}$ by contracting the two
  $-2$ curves above $p_{13(21)^2}$, and the $(-2)$-curve above
  $p_{12}$, the $(-2)$ and the $(-3)$-curves above $p_{13}$ and
  $p_{23}$. We still denote by $\gamma:\widehat{W}\rightarrow W$ the
  contraction, and $\varphi:W\rightarrow X$ the corresponding
  birational transformation. We denote the exceptional lines by $E$,
  $F$, $G$, $H$, above $p_{13(21)^2}$, $p_{12}$, $p_{13}$ and $p_{23}$
  respectively.
\end{dfn}

\begin{rk}
Whenever a point of $x\in X$ is not blown-up in order to get $Y$
(resp. $Z$), we will use the same notation for its proper transform in
$Y$, $\widehat{Y}$, $Z$ or $\widehat{Z}$.
\end{rk}

\begin{thm}\label{thm:main}
\begin{enumerate}
\item The pair $(X',\frac{2}{3}M')$ is a ball quotient orbifold with
  one cusp, where $X'=X\setminus\{p_{13(21)^2}\}$ and $M'=M\cap X'$.
\item The pair $(Y',\frac{3}{4}M'+\frac{3}{4}E)$ is a ball quotient
  orbifold with one cusp, where $Y'=Y\setminus\{p_{12}\}$, and $M'$
  denotes the intersection with $Y'$ of the strict transform of $M$ in
  $Y$.
\item The pair $(Z',\frac{5}{6}M'+\frac{1}{2}E+\frac{5}{6}F)$ is a
  ball quotient orbifold with two cusps, where
  $Z'=Z\setminus\{p_{13},p_{23}\}$, and $M'$ denotes the intersection
  with $Z'$ of the strict transform of $M$ in $Z$.
\item The pair $(W',\frac{1}{2}F+\frac{1}{2}G+\frac{1}{2}H)$ is a ball
  quotient orbifold with one cusp, where $W'=W\setminus M'$, and $M'$
  denotes strict transform of $M$ in $W$.
\end{enumerate}
\end{thm}

\section{Proof of the main result}

The basis of the proof, like in~\cite{derauxklein}, will be a detailed
study of the pairs $(X^{(p)},D^{(p)})$, where $X^{(3)}=X$,
$X^{(4)}=Y$, $X^{(6)}=Z$, $X^{(\infty)}=W$ and the $D^{(p)}$'s are
$\Q$-divisors given by $D^{(3)}=\frac{2}{3}M$,
$D^{(4)}=\frac{3}{4}\widetilde{M}+\frac{3}{4}E$,
$D^{(6)}=\frac{5}{6}\widetilde{M}+\frac{1}{2}E+\frac{5}{6}F$ and
  $D^{(\infty)}=\widetilde{M}+\frac{1}{2}F+\frac{1}{2}G+\frac{1}{2}H$.

\begin{prop}
For each $p$ as above,
\begin{enumerate}\label{prop:main}
\item the pair $(X^{(p)},D^{(p)})$ has at worst log canonical
  singularities;
\item the log-canonical divisor $K_{X^{(p)}}+D^{(p)}$ is ample, i.e.
  the pair $(X^{(p)},D^{(p)})$ is its own canonical model;
\item $c_1^2(X'^{(p)},D'^{(p)})=3c_2(X'^{(p)},D'^{(p)})$, where
  $X'^{(p)}$ denotes the log-terminal locus (obtained from $X^{(p)}$
  by removing the points where the pair is not log-terminal), and
  $D'{(p)}=D^{(p)}\cap X'{(p)}$.∂
\end{enumerate}
\end{prop}
By a theorem of Kobayashi, Nakamura and Sakai~\cite{kns},
Proposition~\ref{prop:main} implies Theorem~\ref{thm:main}.

\begin{figure}
  \hfill
  \begin{subfigure}{0.3\textwidth}
    \includegraphics[width=\textwidth]{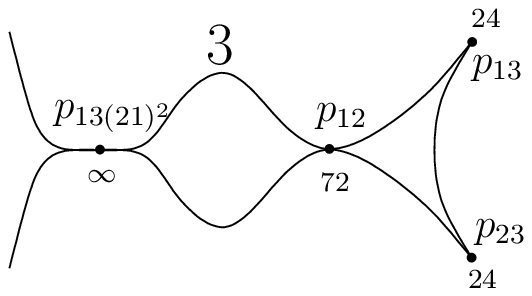}\caption[(a)]{$X$}
  \end{subfigure}\hfill
  \begin{subfigure}{0.3\textwidth}
    \includegraphics[width=\textwidth]{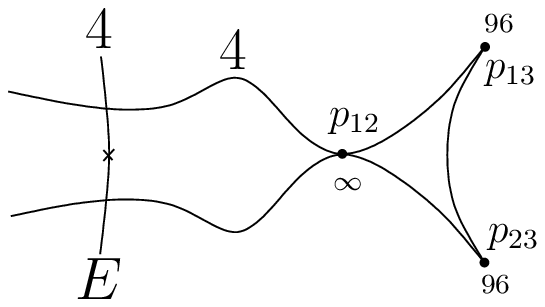}\caption{$Y$}
  \end{subfigure}\hfill\,\\
  \hfill
  \begin{subfigure}{0.3\textwidth}
    \includegraphics[width=\textwidth]{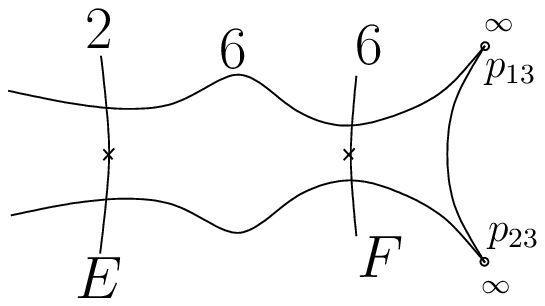}\caption{$Z$}
  \end{subfigure}\hfill
  \begin{subfigure}{0.35\textwidth}
    \includegraphics[width=\textwidth]{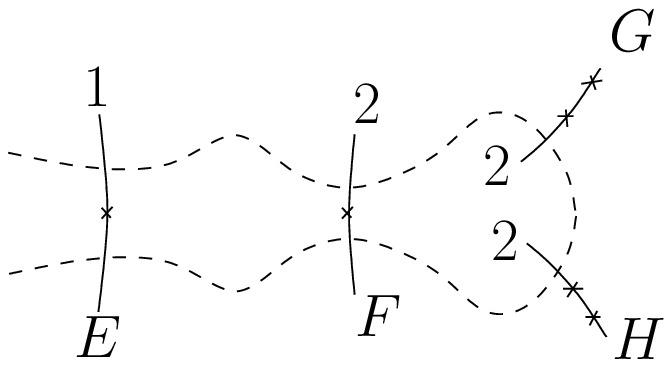}\caption{$W$}
  \end{subfigure}
  \hfill\,
\caption{Schematic picture of the orbifold structure on $X=A/F$, and
  of the relevant birational surfaces $Y$, $Z$ and $W$. We label each
  curve with the relevant orbifold weight (in the case of $W$, the
  dotted curve is removed).}\label{fig:models}
\end{figure}

A schematic picture of the spaces $X$, $Y$, $Z$, $W$ showing the
combinatorics/singularities of the relevant $\Q$-divisors is given in
Figure~\ref{fig:models}. Note that all these spaces map to $X$, and
these maps are isomorphisms over $X\setminus M$, where $X$ (hence
$Y$, $Z$, $W$ as well) has two isolated singularities, of type
$\frac{1}{3}(1,2)$ and $\frac{1}{8}(1,3)$.

\subsection{Log-canonical singularities}

For part~(1) of Proposition~\ref{prop:main}, the only point to
consider is the point $p_{13(21)^2}$, since the others have local
descriptions that were handled in~\cite{derauxklein}. In what follows,
to simplify notation, we write $q=p_{13(21)^2}$.

At the point $q$, we denote by $\widehat{X}$ the minimal resolution of
the pair $(X,\lambda M)$, which is given by
$\pi:\widehat{X}\rightarrow X$, and has exceptional locus a $-1,-2,-2$
chain of projective lines, denoted by $E_1$, $E_2$ and $E_3$ (note
that $E_1$ intersects the proper transform $\widehat{M}$ twice, but
$E_2$ and $E_3$ do not intersect $\widehat{M}$). One checks that
$$
K_{\widehat{X}}=\pi^*K_X+E_1+2E_2+3E_3,
$$
and
$$
\pi^*M=\widehat{M}+2E_1+4E_2+6E_3.
$$
This gives
$$
K_{\widehat{X}}+\lambda\widehat{M}=\pi^*(K_X+\lambda M)+(1-2\lambda)E_1+(2-4\lambda)E_2+(3-6\lambda)E_3,
$$ 
hence the pair $(X,\lambda M)$ is log-canonical at $q$ if and only if
$\lambda \leq 2/3$. For $\lambda=1-1/p$, this means $p\leq 3$. For
$p=3$ the pair is \emph{not} log-terminal at $q$.

Near $p_{12}$, the pair $(X,M)$ is log-canonical for $p\leq 4$, and
log-terminal for $p<4$; at $p_{13}$ and $p_{23}$, it is log-canonical
for $p\leq 6$, log-terminal for $p<6$ (see~\cite{derauxklein} for more
details, where the same type of singularities of the pair occur).

\subsection{Miyaoka-Yau equality} \label{sec:miyaoka}

The formulas for $c_1^2(X^{(p)},D^{(p)})$ are very similar to those
in~\cite{derauxklein}. If we knew that $X$ was a weighted projective
plane, the formulas below would be obtained from those
in~\cite{derauxklein} by replacing 2,3,7 by 1,3,8. We give a slightly
different argument, that relies on the fact that $X=A/F$, where $F$ is
a specific group of order 48. In other words, there is a map
$f:A\rightarrow X$ of degree 48, that ramifies with order 2 along the
union $\cup E_j$ of 12 elliptic curves. 

It follows from the discussion in section~\ref{sec:acrg} that for
every $k$, $E_k\cdot\sum_{j=1}^{12}E_j=24$ (more specifically, see
table~\ref{tab:mirror}). From this, it follows that 
$$
M^2=\frac{1}{48}(2\sum E_j)^2=\frac{1}{48}\cdot 4\cdot 12\cdot 24 = 24.
$$

Note also that $f^*(K_X+\frac{1}{2}M)=K_A$, so 
$$
K_X\cdot M = \frac{1}{48}(K_A\cdot f^*M-\frac{1}{2}(f^*M)^2)=-12,
$$
where we have used the adjunction formula and the fact that $\chi(E_j)=0$.

Finally, note that $(K_X+\frac{1}{2}M)^2=0$ (since $K_A$ is trivial), hence
$$
K_X^2=-K_X\cdot M-\frac{1}{4}M^2=6.
$$

In particular, we get for any $\lambda=1-1/p$, that
$$
(K_X+\lambda M)^2=6(-1+2\lambda)^2=\frac{1}{24}(-12+24\lambda)^2.
$$ 
The last expression is written so as to resemble the formula
in~\cite{derauxklein}.

We now write $D=\lambda M$ on $X$, $\lambda\widetilde{M}+\mu E$ on
$Y$, $\lambda\widetilde{M}+\mu E +\nu F$ on $Z$, $\lambda\widetilde{M}
+ \mu E + \nu F + \sigma G + \sigma H$ on $W$ (recall that the
coefficient of each divisor has the form $1-1/k$, where $k$ is an
integer or $\infty$).

We get for $p=3$,
$$
(K_X+\lambda M)^2 = \frac{1}{24}(-12+24\frac{2}{3})^2 = \frac{2}{3}.
$$
For $p=4$, we take $\lambda=\mu=1-1/4$, and get
$$
(K_Y+\lambda \widetilde{M}+\mu E)^2 = \frac{1}{24}(-12+24\lambda)^2 - \frac{1}{3}(3-6\lambda+\mu)^2 = \frac{21}{16}.
$$
For $p=6$, we take $\lambda=\nu=1-1/6$, $\mu=1-1/2$, and get
$$
(K_Z+\lambda \widetilde{M}+\mu E+\nu F)^2 = 
   \frac{1}{24}(-12+24\lambda)^2 - \frac{1}{3}(3-6\lambda+\mu)^2 - \frac{1}{2}(2-4\lambda+\nu)^2 = \frac{43}{24}.
$$
For $p=\infty$, we take $\lambda=1$, $\mu=1-1=0$, $\nu=\sigma=\tau=1-1/2$, and get
\begin{eqnarray*}
&&(K_W+\lambda \widetilde{M}+\mu E + \nu F +\sigma G +\sigma H)^2 = \\
  && \frac{1}{24}(-12+24\lambda)^2 - \frac{1}{3}(3-6\lambda+\mu)^2 - \frac{1}{2}(2-4\lambda+\nu)^2 
    - 2\cdot \frac{1}{6}(4-6\lambda + \sigma)^2 = \frac{9}{8}.
\end{eqnarray*}

The orbifold Euler characteristics are given by the following. For $p=3$, we get
$$
\chi^{orb}(X,D) = \frac{1}{3} + \frac{1}{8} +\frac{1}{72} + 2\frac{1}{24} + \frac{-4}{3} + 1 = 2/9.
$$
For $p=4$, we get
$$
\chi^{orb}(X,D) = \frac{1}{3} + \frac{1}{8} + 2\frac{1}{4\cdot 4} +
  \frac{1}{3\cdot 4} + 2\frac{1}{96} + \frac{-4}{4} + \frac{-1}{4} + 1 = \frac{7}{16}.
$$
For $p=6$, we get
$$
\chi^{orb}(X,D) = \frac{1}{3} + \frac{1}{8} + 2\frac{1}{2\cdot 6} + 2\frac{1}{6\cdot 6} +
  \frac{1}{3\cdot 2} +\frac{1}{2\cdot 6} + \frac{-4}{6} + \frac{-1}{2} +\frac{-1}{6} + 1 = \frac{43}{72}.
$$

For $p=\infty$, we get
$$
\chi^{orb}(X,D) = \frac{1}{3} + \frac{1}{8} + \frac{1}{3} + \frac{1}{2\cdot 2} +
  2\frac{1}{2\cdot 2} + 2\frac{1}{2\cdot 3} + \frac{-1}{1} + \frac{-1}{2} + 2\frac{-1}{2} + 1 = \frac{3}{8}.
$$

Putting this together, we get that $c_1^2=3c_2$ for all relevant values of $p$.

\subsection{Ampleness}

Our argument relies in part on the following fact, which would be
obvious if we knew $X$ to be a weighted projective plane.
\begin{prop}
  Let $X=\C^2/G=A/F$ be as above. Then $\chi(X)=3$, and $Pic(X)=\Z$.
\end{prop}

\begin{pf}
  The fact that $\chi(X)=3$ follows from the arguments in
  Proposition~\ref{prop:euler}. Indeed, we use the stratification of
  $X$ corresponding to isotropy groups, there are 2 isolated
  singularities, 4 points with non-cyclic reflection stabilizers, 1
  6-punctured projective line, and the open part that has Euler
  characteristic 1. We then have
  $$
  \chi(X)=6-4+1=3.
  $$

  We then use the fact that $X$ is simply connected, because its
  orbifold fundamental group is generated by point stabilizers (this
  follows from a theorem of Armstrong, see~\cite{armstrong}). This
  gives $b_1(X)=0$. Since $X$ has quotient singularities, it satisfies
  Poincar\'e duality (see Theorem~1.13 of~\cite{steenbrink}), hence
  $\chi(X)=3$ gives $b_2(X)=1$.

  The fact that the Picard number is one then follows, see the proof
  of Proposition~4.20 of~\cite{ou}, for instance.
\end{pf}

From this and the analysis in the beginning of
section~\ref{sec:miyaoka}, it follows that $K_X$ is numerically
equivalent to $-\frac{1}{2}M$.

We want to check whether the log-canonical divisors
$K_X+\frac{2}{3}\widetilde{M}$ (case $p=3$),
$K_Y+\frac{3}{4}\widetilde{M}+\frac{3}{4}E$ (case $p=4$),
$K_Z+\frac{5}{6}\widetilde{M}+\frac{1}{2}E+\frac{5}{6}F$ (case $p=6$),
$K_W+\widetilde{M}+\frac{1}{2}F+\frac{1}{2}G+\frac{1}{2}H$ (case
$p=\infty$) are ample.

For the case $p=3$, we simply have $K_X+\frac{2}{3}M\equiv
\frac{1}{6}M$, which is clearly ample (for instance, by the
Nakai-Moishezon it is enough to show that its intersection with $M$ is
$>0$, but $M^2=24>0$).

For $p=4$, we have
$$
K_Y+\lambda\widetilde{M}+\mu E = \varphi^*(K_X+\lambda M)+(3-6\lambda+\mu)E,
$$ 
where $\varphi$ is as in Definition~\ref{def:Y}. Since
$\varphi^*M\equiv \widetilde{M}+6E$, the right hand side is linearly
equivalent to
$$
(\lambda-\frac{1}{2})(\widetilde{M}+6E)+(3-6\lambda+\mu)E=(\lambda-\frac{1}{2})\widetilde{M}+\mu
E.
$$ 
We check that the latter divisor is ample by the Nakai-Moishezon
criterion. As explained in section~3.3 of~\cite{derauxklein}, since
$Pic(X)=\Z$, it is enough to check that its intersection with
$\widetilde{M}$ and with $E$ is $>0$.

Now we recall that $E^2=-1/3$, and compute
\begin{eqnarray*}
&&(\varphi^*(\lambda -\frac{1}{2})M + (3-6\lambda+\mu)E)\cdot \widetilde{M}\\
&&=(\varphi^*(\lambda -\frac{1}{2})M + (3-6\lambda+\mu)E)\cdot (\varphi^*M-6E)\\
&&=(\lambda-\frac{1}{2})M^2-6(3-6\lambda+\mu)E^2=\frac{9}{2}>0,
\end{eqnarray*}
and
$$
(\varphi^*(\lambda -\frac{1}{2})M + (3-6\lambda+\mu)E)\cdot E = (3-6\lambda+\mu)E^2=\frac{1}{4}>0.
$$

The cases $p=6$, $p=\infty$ are similar, simply with slightly longer
computations. The basis of the computation is
$$ 
K_Z+\lambda\widetilde{M}+\lambda E + \mu F \equiv
(\lambda-\frac{1}{2})\varphi^*M + (3-6\lambda+\mu)E +
(2-4\lambda+\nu)F,
$$
and
{\footnotesize
\begin{eqnarray*}
&&K_Z+\lambda\widetilde{M}+\lambda E + \mu F +\sigma G + \sigma H \\
&&\equiv
(\lambda-\frac{1}{2})\varphi^*M + (3-6\lambda+\mu)E +
(2-4\lambda+\nu)F + (4-6\lambda+\sigma)(G+H).
\end{eqnarray*}
}
Also, for $\varphi:Z\rightarrow X$, we have 
$$
\varphi^*M=\widetilde{M}+6E+4F,
$$
and for $\varphi:W\rightarrow X$, we have 
$$
\varphi^*M=\widetilde{M}+6E+4F+6(G+H).
$$

\subsection{Identifying the groups}

In this section we briefly explain why the holonomy group of the
complex hyperbolic structures constructed by uniformization (using the
Kobayashi-Nakamura-Sakai version of equality case in the Miyaoka-Yau
inequality) is isomorphic to the relevant sporadic triangle groups.

\begin{thm}
  Let $\Gamma_p$ be the group obtained from the statement of
  Theorem~\ref{thm:main} for $p=3$, $4$, or $6$. Then $\Gamma_p$ is
  conjugate to the triangle sporadic group $\Sc(p,\sigma_1)$.
\end{thm}

\begin{pf}
  This follows from the description of orbifold fundamental group
  $\Gamma_2$, which is generated by three complex reflections $R_j$,
  $j=1,2,3$ of order 2, such that:
  \begin{itemize}
    \item $\br(R_1,R_2)=4$, $\br(R_2,R_3)=3$,
      $\br(R_3,R_1)=3$.
    \item $R_1R_2R_3$ (has linear part which) is regular elliptic of
      order 8.
  \end{itemize}
  Given how the orbifold structure with holonomy $\Gamma_p$ is
  constructed, these same properties will hold, with complex
  reflections of order $p$ instead of order 2, except that in the
  cases $p>2$, the isometry $R_1R_2R_3$ is regular elliptic of order 8
  (the analogue of taking the linear part is then simply to view it as
  an element of the stabilizer of its fixed point, which is isomorphic
  to $U(2)$).

  The result then follows from Proposition~\ref{prop:unique}, stated
  and proved below.
\end{pf}

\begin{prop}\label{prop:unique}
  Let $\Gamma$ be a lattice generated by three complex reflections
  $R_j$, $j=1,2,3$ such that ${\textrm br}(R_1,R_2)=4$, ${\textrm
    br}(R_2,R_3)=3$, ${\textrm br}(R_3,R_1)=3$ and $R_1R_2R_3$ is
  regular elliptic of order 8. Then $\Gamma$ is conjugate to the
  Thompson group $\T(p,{\bf E_1})$, which is isomorphic to the
  sporadic triangle group $\Sc(p,\sigma_1)$.
\end{prop}
\begin{pf}
  First, the fact that $\T(p,{\bf E_1})$ and $\Sc(p,\sigma_1)$ are
  conjugate follows from a suitable change of generators, in the same
  vein as in~\cite{kamiyaparkerthompson}; the details are given in
  Proposition~7.1 of~\cite{thealgo}.

  Any group as above must be conjugate to $\Sc(p,\sigma_1)$ or
  $\Sc(p,\bar\sigma_1)$, but the last group is not discrete if $p=3$ or
  $6$ (see section Section~9.4 of~\cite{dpp1}).

  One checks that $\Sc(4,\bar\sigma_1)$ is not discrete either, for
  instance by showing that $M=R_3R_1R_2J$ is regular elliptic but has
  infinite order (here $R_1$, $R_2$ and $R_3$ stand for the standard
  generators of $\Sc(4,\bar\sigma_1)$, and $J$ stands for the regular
  elliptic element of order 3 that conjugates $R_j$ into $R_{j+1}$).
  Indeed, one checks that $\tr(
  M)=(\sqrt{3}+i)(i-(1+i)\sqrt{2})/2:=\tau$, and
  $$
  |\tau|^4-8\mathfrak{Re}(\tau^3)+18|\tau|^2-27=88-64\sqrt{2}<0,
  $$ 
  so $M$ is regular elliptic (see~\cite{goldman}, section~6.2.3).

  The characteristic polynomial of $M$ is equal to
  $\lambda^3-\tau\lambda^2+\bar\tau\lambda-1$, and one verifies that
  only one of its roots is a root of unity (namely
  $-(i+\sqrt{3})/2$). Indeed, the other two roots have a minimal
  polynomial of degree 16, that is not cyclotomic.

  Note that in $\Sc(4,\bar\sigma_1)$, the element $R_3R_1R_2J$ is
  loxodromic (and indeed, we know that this is a lattice,
  see~\cite{thealgo}).
\end{pf}

\begin{rk}
  The argument we just gave provides a short proof that
  $\Sc(p,\sigma_1)$ is indeed a lattice for $p=3,4,6$, a fact which
  was proved using heavy computer power in~\cite{thealgo}.
\end{rk}

For $p=\infty$, Proposition~\ref{prop:unique} has the following
analogue (recall that unipotent elements are isometries whose matrix
representative has a single eigenvalue of multiplicity 3), which can
be proved with very similar methods as in~\cite{dpp1}. We omit the 
details because the corresponding group turns out to be arithmetic.
\begin{prop}\label{prop:unique2}
  Let $\Gamma$ be a lattice generated by unipotent elements $R_j$,
  $j=1,2,3$ such that ${\textrm br}(R_1,R_2)=4$, ${\textrm
    br}(R_2,R_3)=3$, ${\textrm br}(R_3,R_1)=3$ and $R_1R_2R_3$ is
  regular elliptic of order 8. Then $\Gamma$ is conjugate to Thompson
  group $\T(\infty,{\bf E_1})$, which is isomorphic to the sporadic
  triangle group $\Sc(\infty,\sigma_1)$.
\end{prop}
The arithmeticity of the group $\Sc(\infty,\sigma_1)$ is fairly
obvious from the description given in section~\ref{sec:sporadic},
where we give a generating set with entries in $\Z+ i\sqrt{2}\Z$. From
this it follows that the adjoint trace field is $\mathbb{Q}$, hence
the group is indeed arithmetic (see~\cite{thealgo}, for instance).


\begin{thebibliography}{10}

\bibitem{armstrong}
M.~A. Armstrong.
\newblock The fundamental group of the orbit space of a discontinuous group.
\newblock {\em Proc. Cambridge Philos. Soc.}, 64(2):299--301, 1968.

\bibitem{bannai}
E.. Bannai.
\newblock Fundamental groups of the spaces of regular orbits of the finite
  unitary reflection groups of dimension 2.
\newblock {\em J. Math. Soc. Japan}, 28:447--454, 1976.

\bibitem{bhh}
G.~Barthel, F.~Hirzebruch, and T.~H{\"o}fer.
\newblock {\em Geradenkonfigurationen und {A}lgebraische {F}l\"achen}.
\newblock Aspects of Mathematics, D4. Friedr. Vieweg \& Sohn, Braunschweig,
  1987.

\bibitem{bernsteinSchwarzman}
J.~Bernstein and O.~Schwarzman.
\newblock Chevalley's theorem for the complex crystallographic groups.
\newblock {\em J. Nonlinear Math. Phys.}, 13(1-4):323--351, 2006.

\bibitem{birkenhakelange}
C.~Birkenhake and H.~Lange.
\newblock {\em {Complex Abelian Varieties}}.
\newblock Springer, 2004.

\bibitem{chengreenberg}
S.~S. Chen and L.~Greenberg.
\newblock Hyperbolic spaces.
\newblock {Contribut. to Analysis, Collect. of Papers dedicated to Lipman Bers,
  49-87 (1974)}, 1974.

\bibitem{chl}
W.~{Couwenberg}, G.~{Heckman}, and E.~{Looijenga}.
\newblock Geometric structures on the complement of a projective arrangement.
\newblock {\em {Publ. Math., Inst. Hautes \'Etud. Sci.}}, 101:69--161, 2005.

\bibitem{delignemostow}
P.~Deligne and G.~D. Mostow.
\newblock Monodromy of hypergeometric functions and non-lattice integral
  monodromy.
\newblock {\em Inst. Hautes {\'E}tudes Sci. Publ. Math.}, 63:5--89, 1986.

\bibitem{delignemostowbook}
P.~Deligne and G.~D. Mostow.
\newblock {\em Commensurabilities among lattices in {PU}(1,n)}, volume 132 of
  {\em Annals of Mathematics Studies}.
\newblock Princeton Univ. Press, Princeton, 1993.

\bibitem{derauxklein}
M.~Deraux.
\newblock {Non-arithmetic lattices and the Klein quartic}.
\newblock to appear in J. reine Angew. Math., arXiv:1605.03846,.

\bibitem{thealgo}
M.~Deraux, J.~R. Parker, and J.~Paupert.
\newblock On commensurability classes of non-arithmetic complex hyperbolic
  lattices.
\newblock Preprint, arXiv:1611.00330.

\bibitem{dpp1}
M.~Deraux, J.~R. Parker, and J.~Paupert.
\newblock Census for the complex hyperbolic sporadic triangle groups.
\newblock {\em Experiment. {M}ath.}, 20:467--486, 2011.

\bibitem{dicerbostover}
L.~F. Di~Cerbo and M.~Stover.
\newblock Bielliptic ball quotient compactifications and lattices in {PU}(2,1)
  with finitely generated commutator subgroup.
\newblock {\em Ann. Inst. Fourier (Grenoble)}, 67(1):315--328, 2017.

\bibitem{dolgachev}
I.~V. Dolgachev.
\newblock Reflection groups in algebraic geometry.
\newblock {\em {Bull. Am. Math. Soc., New Ser.}}, 45(1):1--60, 2008.

\bibitem{goldman}
W.~M. Goldman.
\newblock {\em Complex {H}yperbolic {G}eometry}.
\newblock Oxford Mathematical Monographs. Oxford University Press, 1999.

\bibitem{hirzabelian}
F.~Hirzebruch.
\newblock Chern numbers of algebraic surfaces: an example.
\newblock {\em Math. Ann.}, 266(3):351--356, 1984.

\bibitem{kamiyaparkerthompson}
S.~Kamiya, J.R. Parker, and J.M. Thompson.
\newblock Notes on complex hyperbolic triangle groups.
\newblock {\em Conformal Geometry and Dynamics}, 14:202--218, 2010.

\bibitem{kanekotokunagayoshida}
J.~Kaneko, S.~Tokunaga, and M.~Yoshida.
\newblock {Complex crystallographic groups II}.
\newblock {\em J. Math. Soc. Japan}, 34:595--605, 1982.

\bibitem{kns}
R.~Kobayashi, S.~Nakamura, and F.~Sakai.
\newblock A numerical characterization of ball quotients for normal surfaces
  with branch loci.
\newblock {\em Proc. Japan Acad., Ser. A}, 65(7):238--241, 1989.

\bibitem{livne}
R.~A. Livn{\'e}.
\newblock {\em On certain covers of the universal elliptic curve}.
\newblock PhD thesis, Harvard University, 1981.

\bibitem{mostowihes}
G.~D. Mostow.
\newblock Generalized {P}icard lattices arising from half-integral conditions.
\newblock {\em Inst. Hautes {\'E}tudes Sci. Publ. Math.}, 63:91--106, 1986.

\bibitem{ou}
W.~Ou.
\newblock Lagrangian fibrations on symplectic fourfolds.
\newblock To appear in J. Reine Angew. Math., \verb|arXiv:1411.1377|.

\bibitem{popov}
V.~L. Popov.
\newblock Discrete complex reflection groups.
\newblock Commun. Math. Inst., Rijksuniv. Utrecht 15, 1982.

\bibitem{roulleauurzua}
X.~Roulleau and G.~Urz\'ua.
\newblock Chern slopes of simply connected complex surfaces of general type are
  dense in [2,3].
\newblock {\em Ann. Math. (2)}, 182(1):287--306, 2015.

\bibitem{shephardtodd}
G.~C. Shephard and J.~A. Todd.
\newblock Finite unitary reflection groups.
\newblock {\em Canadian J. Math.}, 6:274--304, 1954.

\bibitem{springer}
T.~A. Springer.
\newblock {\em Invariant {T}heory}, volume 585 of {\em Lecture Notes in
  Mathematics}.
\newblock Springer-Verlag, Berlin-Heidelberg-New York, 1977.

\bibitem{steenbrink}
J.~H.~M. Steenbrink.
\newblock Mixed {H}odge structure on the vanishing cohomology.
\newblock Real and compl. Singul., Proc. Nordic Summer Sch., Symp. Math., Oslo
  1976, 1977.

\bibitem{stover}
M.~Stover.
\newblock Cusp and $b_1$ growth for ball quotients and maps onto $\mathbb{Z}$
  with finitely generated kernel.
\newblock Preprint, arXiv:1506.06126.

\bibitem{thompson}
J.~M. Thompson.
\newblock {\em Complex Hyperbolic Triangle Groups}.
\newblock PhD thesis, Durham University, 2010.

\end{thebibliography}
\end{document}